\documentclass[english]{article}
\usepackage[T1]{fontenc}
\usepackage[latin1]{inputenc}
\usepackage{textcomp}
\usepackage{amsmath}
\usepackage{amssymb}
\usepackage{stmaryrd}
\usepackage{babel}
\begin{document}

\author{Anne Bertrand Mathis}

\author{Université de Poitiers (France)}

\author{anne.bertrand@math.univ-poitiers.fr }

\title{Democratic sequences}

\title{Algorithmic construction of a sequence taking values in a coutable
alphabet with frequency preset.}

\author{Anne Bertrand Mathis}

\author{Université de Poitiers (France)}

\author{anne.bertrand@math.univ-poitiers.fr }
\maketitle
\begin{abstract}
Let $A=\left\{ a_{n};n\geq1\right\} $ be an enumerable alphabet,
and let $\left\{ \lambda_{n};\thinspace n\geq1\right\} $ be a sequence
of positive numbers such that $\sum_{n\geq1}\lambda_{n}=1.$ We explain
an algorithmic construction of a sequence $\left(u_{n}\right)_{n\geq1}$
on the alphabet $A$ in which each letter $a_{k}$ appears with the
frequency $\lambda_{k}=\lambda\left(a_{k}\right).$

Etant donné un alphabet fini ou dénombrable $A=\left\{ a_{n};n\geq1\right\} $
et une famille $\left\{ \lambda_{n};n\geq1\right\} $ de nombres positifs
de somme $1$ nous construisons de façon simple une suite $\left(u_{n}\right)_{n\geq1}$
à valeurs dans $A$ dans laquelle chaque lettre $a_{k}$ apparaît
avec la fréquence $\lambda_{k}=\lambda\left(a_{k}\right)$.


Classification AMS : 11K99
\end{abstract}

\section{Presentation of the method.}

\textbf{Method}: Let $A$ be a finite or countable alphabet and $\left(\lambda\left(a\right)\right)_{a\in A}$
a family of positive numbers with sum $1.$ We order the alphabet
$A$: $A=\left\{ a_{1},a_{2},...a_{p},...\right\} $ so that if $h$
is less than $k$ then $\lambda\left(a_{h}\right)\leq\lambda\left(a_{k}\right)$
; $p$ is called the index of $a_{p}.$ 

We arbitrarily choose the start $u_{1}...u_{k}$ of the sequence $\left(u_{n}\right)_{n\geq1}$
that we want to construct. It can be a single letter $u_{1}.$ Let
us assume the first $M$ terms $u_{1},...,u_{M}$ already constructed.
We call $M-$frequency of a letter $a\in A$ the number $\lambda_{M}\left(a\right)=\frac{1}{M}\sum_{u_{n}=a,n\leq M}1$.
We call \textit{$M-$deficit }of $a$ the number $D_{M}\left(a\right)=\lambda\left(a\right)-\lambda_{M}\left(a\right)$
when positive. If it is strictly positive we say that $a$ is \textit{late
at step $M$ }and if $D_{M}\left(a\right)>$$D_{M}\left(b\right)$
we say that at step $M$ $a$ is more late than $b$.

If $\lambda\left(a\right)-\lambda_{M}\left(a\right)$ is negative
we say that $a$ is \textit{ahead }in step $M$. We then call $M-$\textit{excess}
de $a$ at step $M$ the strictly positive number 

$E_{M}\left(a\right)=\mid\lambda\left(a\right)-\lambda_{M}\left(a\right)\mid=\lambda_{M}\left(a\right)-\lambda\left(a\right).$ 

Once $u_{1},...,u_{M}$ determined we choose for $u_{M+1}$ that of
the letters wich is the most late or, if several present the same
delay, the one wich has the lowest index (we can also choose at random
one of the most late letters, they are in finite numbers). If none
of the letters are late we choose $a_{1}.$ 

\textbf{Théorème} 1. \textit{Let $A$ be a finite or countable alphabet,
$\left(\lambda\left(a\right)\right)_{a\in A}$ a family of positive
numbers with sum $1,$ and let $\left(u_{n}\right)_{n\geq1}$ be a
constructed sequence accorded to the method above.}

\textit{Then any letter $a$ of the alphabet $A$ appears in the sequence
$\left(u_{n}\right)_{n\geq1}$ with the frequency $\lambda\left(a\right)$
:}

\[
\forall a\in A,\ \ \lim_{N\rightarrow\infty}\frac{1}{N}\sum_{_{\begin{array}{c}
n\leq N\\
u_{n}=a
\end{array}}}1=\lambda\left(a\right)\text{  }(1).
\]

This method is inspired by the work of Ville $\left[3\right]$.

The alphabet does not required to be countable, the condition $\sum_{a\in A}\left(\lambda\left(a\right)\right)=1$
implies that only a countable amount of letters have a stricly positive
frequency.

Furthermore we can also start the sequence \textit{$\left(u_{n}\right)_{n\geq1}$}
with a Joker $J$ of zero frequency wich will only appears once and
has frequency zero, this is what we do if we want the first terms
of the sequence $\left(u_{n}\right)_{n\geq1}$ are not chosen arbitrarily
(this will be the case in the study of numeration sequences). 

The proof of the theorem is simple when the number of letters is finite
but if the letters are in countable quantity the situation becomes
complicated. Let us state some lemmas.

\textbf{Lemma 1.} \textit{$a)$ For any integer entier $M$ the sum
of the deficits of the letters of $A$ is equal to the sum of the
excesses of the letters of $A$.}

\textit{$b)$ If there exists a letter $b$ in advance with $E_{M}\left(b\right)$
strictly positive then there exists a late letter $a$ with $D_{M}\left(a\right)$
strictly positive and vice versa. }

\textit{$c)$ The algebraic sum of the deficits and excesses of the
letters of $A$ is equal to zero:}

\[
\sum_{a\in A}D_{M}\left(a\right)-\sum_{a\in A}E_{M}\left(a\right)=\sum_{a\in A}\lambda\left(a\right)-\lambda_{M}\left(a\right)=0.
\]

\textit{Proof of assertion $c)$.} Immediate from the definitions
of deficits and excesses: in fact the numbers of occurrences of the
letters of the alphabet between the step $1$ and $M$ is equal to
$M$ therefore $\sum_{a\in A}\frac{\lambda_{M}\left(a\right)}{M}=1=\sum_{a\in A}\left(\lambda\left(a\right)\right)$
and \textit{$\sum_{a\in A}D_{M}\left(a\right)-\sum_{a\in A}E_{M}\left(a\right)=\sum_{a\in A}(\lambda\left(a\right)-\lambda_{M}\left(a\right))=0.$ }

Let us call the \textit{upper frequency of a letter} $b$ the number
$\overline{\lambda}\left(b\right)=\limsup_{N\rightarrow\infty}\frac{1}{N}\sum1$
where the sum is taken on on the integers $n\leq N$ such that $u_{n}=b$.
We also define the \textit{lower frequency} of a letter. If the sequence
$\left(u_{n}\right)_{n\geq1}$constructed according to the method
indicated does not satisfy condition $\left(1\right)$ it is because
there exists one or more letters $b$ belonging to $A$ wich do not
admit $\lambda\left(b\right)$ for frequency. Therefore there exists
ar least one letter $b$ whose upper frequency $\overline{\lambda}\left(b\right)$
is strictly greater than $\lambda\left(b\right)$ or a letter $b$
whoses lower frequency is strictly less than $\underline{\lambda}$$\left(b\right)$.
As the number of letters in infinite we are not very sure that the
two cases occur simultaneously.

Remark: perhaps the assertions in lemma $1$should only be valid after
a certain rank due to the selection of the first terms.

\textbf{Proposition 1.}\textit{ For any letter of the alphabet $A$,
$\overline{\lambda}\left(b\right)\leq\lambda\left(b\right)$}\textbf{ }

Suppose that there exists a letter $b$ wich satisfies \textbf{$\overline{\lambda}\left(b\right)>\lambda\left(b\right)$.}

Let us then write $\overline{\lambda}\left(b\right)=\lambda\left(b\right)+f_{b}$
with of course $f_{b}>0$. There exists a non-empty subset $B$ of
$\mathbb{N}$ such that the frequency of $b$ when $N$ tends to infinity
while remaining in $B$ is equal to $\overline{\lambda}\left(b\right)$
and we can find an integer $N_{0}$ such that:

\[
\forall N>N_{0},N\in B,\ \ \lim_{N\rightarrow\infty}\frac{1}{N}\sum_{_{\begin{array}{c}
n\leq N\\
u_{n}=b
\end{array}}}1>\bar{\lambda}\left(b\right)-\frac{f_{b}}{2}=\lambda\left(b\right)+\frac{f_{b}}{2}\text{  }(2).
\]

As $b$ appears an infinite number of times in the sequence \textit{$\left(u_{n}\right)_{n\geq1}$
}there exists an integer $n_{0}$ verifying $\frac{1}{n_{0}-1}<\frac{f_{b}}{4}$
such that $u_{n_{0}}=b.$ 

Let $M>\max\left(N_{o},n_{0}\right)$ be an integer belonging to $B$.
We therefore have, according to \textit{$\left(2\right)$, }

\[
\frac{1}{M}\sum_{_{\begin{array}{c}
n\leq M\\
u_{n}=b
\end{array}}}1>\lambda\left(b\right)+\frac{f_{b}}{2}\text{  }.
\]

Let's go back to the las index $m$ less than or equal to $\text{M }$
such that $u_{m}=b.$ This index $m$ is greater than or equal to
$n_{0}.$ The number of occurences of $b$ among the fist $m$ letters
of the sequence $\left(u_{n}\right)_{n\geq1}$ is equal to the number
of occurrences of $b$ among the first $M$ letters. As more $\frac{M}{m}$
is greater or equal to $1$ it comes:

\[
\frac{1}{m}\sum_{_{\begin{array}{c}
n\leq m\\
u_{n}=b
\end{array}}}1=\text{  }\frac{1}{m}\sum_{_{\begin{array}{c}
u_{n}=b\\
n\leq M
\end{array}}}1=\frac{M}{m}\left(\frac{1}{M}\sum_{_{\begin{array}{c}
n\leq M\\
u_{n}=b
\end{array}}}1\text{  }\right)>\lambda\left(b\right)+\frac{f_{b}}{2}\text{  }\text{  }
\]

Let's look at what happens at step $m-1.$ Let us show that at this
step $b$ was already in advance. The number of occurrences of $b$
when $n\leq m-1$ is equal to the number of occurrences of $b$ betweens
steps $1$ and $m$, minus $1$. Therefore 
\[
\frac{1}{m-1}\sum_{_{\begin{array}{c}
n\leq m-1\\
u_{n}=b
\end{array}}}1=\frac{1}{m-1}\left(\left(\sum_{_{\begin{array}{c}
n\leq m\\
u_{n}=b
\end{array}}}1\right)-1\right)\geq\frac{1}{m}\left(\sum_{_{\begin{array}{c}
n\leq m\\
u_{n}=b
\end{array}}}1\right)-\frac{1}{m-1}\geq
\]

$\geq\lambda\left(b\right)+\frac{f_{b}}{2}-\frac{1}{m-1}.$

As $m$ is greater than or equal to $n_{0}+1,$ $\frac{1}{m-1}$ is
smaller than $\frac{f}{4}$ and at the step $m-1$ the letter $b$
admits an excess of at least $\frac{f_{b}}{2}-\frac{f_{b}}{4}$$=\frac{f_{b}}{4}$. 

The sum of deficits ans excesses being zero at each step according
to Lemma $1,$ it exists at step $m-1$ a letter $a$ such that $D_{m-1}\left(a\right)$
is strictly positive. So the letter $b$ could not be chosen at this
step, this it impossible so no letter can have an upper density $\overline{\lambda}\left(b\right)$
greater than \textit{$\lambda\left(b\right)$$\boxempty$. }

It remains to show that no letter admits a lower frequency strictly
lower than its assigned frequency.

\textbf{Proposition 2 }\textit{All letter of the alphabet appears
at least once in the sequence $\left(u_{n}\right)_{n\geq1}.$}

We will assume in the following that the letters of the alphabet are
ordered by decreasing frequency: $i<j\Rightarrow\lambda\left(a_{i}\right)\leq\lambda\left(a_{j}\right)$.

\textbf{Lemma 2.}\textit{ The deficit at step $M$ of a letter $a$
is always strictly smaller than $\lambda\left(a\right).$ Therefore
the sum of the deficits of all the letters is at most $1.$ }

\textit{Proof of lemma} $2$. At step $M$ the deficit $D_{M}\left(a\right)$
of the letter $a$ is equal to $\lambda\left(a\right)-\lambda_{M}\left(a\right)$
when his term is positive; as $\lambda_{M}\left(a\right)$ is always
positive or zero, this term is therefore always less than or equal
to $\lambda\left(a\right)$ and so the sum of the deficit is at most
$1.$

\textbf{Lemme 3.}\textit{ Suppose that at a step $M$ the deficit
of the letter $b$ is greater or equal to $g$ ; then no letter $a$
with $\lambda\left(a\right)<g$ can proceed to this step $M$.}

\textit{Prof of lemma $3$. }According to lemma\textit{ $2$} the
deficit of $a$ would then be less than or equal to $\lambda\left(a\right),$
wich is strictly less than the deficit of $b$, $a$ cannot pass.

\textit{Proof of Proposition} $2.$

Suppose there is a letter $a_{i}$ wich never pass. It is this letter
that we will consider. It is assumed that $\lambda\left(a_{1}\right)\geq\lambda\left(a_{2}\right)\geq...\geq\lambda\left(a_{i}\right)\geq\lambda\left(a_{i+1}\right)....$.
.

Suppose that if several letters admit the same frequency as $a_{i}$
their indices are $i-k,...i-1,i$. None of the letters $a_{j}$ such
that $\lambda\left(a_{j}\right)<\lambda\left(a_{i}\right)$ can pass,
therefore none of the letters $a_{i+1},a_{i+2},...$ will pass because
$\lambda\left(a_{i+1}\right)<\lambda\left(a_{i}\right)$ and so tight
now. Therefore all their excess are equal to zero and their deficits
are all strictly lower than their assigned frequency according to
Lemma $2$.

Let $N_{0}$ be such that if $N\geq N_{0}$ the sum of the excesses
of the letters $a_{1},a_{2},...,a_{i-k},...a_{i}$ is less than $\frac{1}{4}\lambda\left(a_{i}\right)$.
These letters are finite in number so proposition $1$ ensures the
existence of $N_{0}$.

Lemma $1a$ show that as long as $a_{i},a_{i+1}...$ have never passed
the sum of the excess of letters $a_{1},a_{2},...,a_{i-k},...,a_{i-1}$
is equal to the sum of their deficits increased by the sum of the
deficits of the letters $a_{i+1},a_{i+2},...$ . The sum of these
deficits is therefore less than $\frac{1}{4}\lambda\left(a_{i}\right)$,
so after step $N_{0}$ the deficits of the letters $a_{1},a_{2},...,a_{i-k},...,a_{i-1}$
cannot be greater than $\frac{1}{4}\lambda\left(a_{i}\right)$. 

For $N>N_{0}$ the sum of excesses is less than $\frac{1}{4}\lambda\left(a_{i}\right)$,
so no letter can have a deficit greater than $\frac{1}{4}\lambda\left(a_{i}\right)$
at step $N.$ But if $a_{i}$ never pass his deficit tend to $\lambda\left(a_{i}\right)$
as $N$ goes to infinity. This is impossible so $a_{i}$ appears in
the sequence.

The following lemma is obvious:

\textbf{Lemma $4$} \textit{Whe the letter $a_{i}$ passes for the
first time to stage $N_{i}$ all the other letters have a deficit
less or equal to $\lambda\left(a_{i}\right).$}

\textbf{Proposition} \textit{$3.$ No letter $a$ can admit a lower
frequency strictly smaller than its assigned frequency : $\forall a\in A$,
$\underline{\lambda}\left(a\right)$ $\geq\lambda\left(a\right).$}

\textbf{Lemma $5$. }\textit{For any integer $M$ and any letter $b,$
}$\mid D_{M+1}\left(b\right)-D_{M}\left(b\right)\mid\leq\frac{1}{M+1}$
\textit{each time these two quantities are defined.}

\textit{Proof of lemma }$5$. When we add a $b$ to step $M$ the\textbf{
}deficit of $b$ decreases of $\frac{(number\thinspace of\thinspace b\thinspace between\thinspace1\thinspace and\thinspace M)+1}{M+1}-\frac{number\thinspace of\thinspace b\thinspace beyween\thinspace1\thinspace and\thinspace M}{M}=\frac{M-number\thinspace of\thinspace b\thinspace between\thinspace1\thinspace and\thinspace M}{M(M+1)}$
wich is smaller than $\frac{1}{M+1}$. When we add a letter different
from $b$ the variation of the deficit is less than $\frac{number~of~b}{M}-\frac{number~of~b}{M+1}=\frac{number~of~b}{M(M+1)}$
who is also less than $\frac{1}{M+1}$.

\textit{Proof of proposition $3$}. Let $B$ be the set of numbers
$\gamma$ de $\left]0,1\right]$ such that there exists $a_{\gamma}$
whose lower frequency is strictly lower than $\lambda\left(a_{\gamma}\right)$
: $\underline{\lambda}\left(a_{\gamma}\right)$ $=\lambda\left(a_{\gamma}\right)-\gamma$.
As the sum of the frequencies is equal to $1$, the number of letters
whose frequency exceed a fixed number $d$ is finite and $B$ cannot
admit an accumulation point except perhaps zero $0$ : if $B$ is
non-empty $B$ admits a maximum element $c$ so there exist at least
one letter $a_{c}$ such that $\underline{\lambda}\left(a_{c}\right)$
$=\lambda\left(a_{c}\right)-c$ and perhaps a finite number of such
letters.

Suppose there is only one letter, $a_{c}$, such that $\underline{\lambda}(a_{c})$$=\lambda\left(a_{c}\right)-c$.
Let $P$ be the (finite) set of letters $a$ such that $\lambda\left(a\right)-$$\underline{\lambda}(a)$$\in\left]0,c\right[$
and let $b$ be the maximum of the quantities $\underline{\lambda}(a)$
for $a\in P$. So $b<c$. Let $Q$ be the set of letters whose frequency
is greater or equal to $b$ and such that $\underline{\lambda}(a)=\lambda\left(a\right)$(and
is also equal to $\bar{\lambda}(a)$, see Proposition $1$); Q is
finite because the sum of the frequencies is $1.$ 

There exists an integer $N_{0}$ such that for all $N>N_{0}$, for
any letter $w$ in $P\cup Q,$ the deficit $D_{N}\left(w\right)$
of the letter $w$ is less than $b+\frac{c-b}{2}$, those of $P$
because their defects have an upper limit at most $b$ and those of
the finite set $Q$ because if $N$ is large enough they are closed
of their assigned frequency ($\underline{\lambda}(a)=\lambda\left(a\right)=\bar{\lambda}(a)$
if $a\in Q$ and $Q$ is finite, use proposition $1$).

Let us place ourselve at a step $N_{i}$ with $N_{i}>N_{0}$ where
the letter $a_{i}$ with $\lambda\left(a_{i}\right)<\frac{b}{2}$
passes for the first time, and such that the\textbf{ }deficit of $a_{c}$
exceed the value $c-\frac{c-b}{20}$ at a step between steps $N_{i}$
et $N_{i+1}$ (such $N_{i}$ exist because the sequence $N_{i}$ tends
toward infinity). The letter $a_{c}$ admits at step $N_{i}$ a deficit
lower than $\frac{b}{2}$ accrding to lemma $4$.The deficit of $a_{c}$
must therefore increase between $N_{i}$ et $N_{i+1},$ very gradually
due to lemme $5$. 

Suppose that at step $M\in\left[N_{i},N_{i+1}\right]$ the deficit
of $a_{c}$ is between $c-\frac{c-b}{10}$ and $c-\frac{c-b}{20}$
; then the letter $a_{c}$ must pass because no other lette has a
greater deficit (lemma $3)$, and it passes until the deficit of $a_{c},$
wich decrease, returns under $c-\frac{c-b}{20}$ ; this deficit will
therefore never be able to rise to exceed $c-\frac{c-b}{20}$. So
the hypothesis is impossible, there cannot be a single letter such
that the deficit of this letter $a_{c}$ admits upper limit $c$;
so the letter $a_{c}$ cannot have a lower frequency $\underline{\lambda}\left(a_{c}\right)$
$=\lambda\left(a_{c}\right)-c$.

Now look at the case where the maximum $c$ of $\lambda\left(a\right)-\underline{\lambda}\left(a\right)$
is reached for several letters, and where all the deficits of the
other letters are smaller than a common value $b<c$. The sets $B,$$P,$$Q$
and the numbers $b<$ $c$ are defined as above and let $R$ designates
the finite set of letters $a$ of the alphabet such that $\underline{\lambda}\left(a\right)=\lambda\left(a\right)-c$
of $R$. 

We also choose for $N_{i}$ an integer as lardge as necessary wich
represents the first passage of a letter $a_{i}$ of assigned frequency
less than the common value $b$ and such that at least one of the
letters of $R$ such that $\underline{\lambda}\left(a\right)=\lambda\left(a\right)-c$
has a deficit exceeding $c-\frac{c-b}{20}$ at a step between steps
$N_{i}$ and $N_{i+1}$(such $N_{i}$ must exist). Letters of $P$
and $Q$ always have deficits less than $b+\frac{c-b}{2}$ for all
$N>N_{i}$. At least one of the letters of $R$, say $w,$ see his
deficit starting from very low ( lemma $4$) and slowly rising to
$c-\frac{c-b}{20}$; it passes if no letter has a greater deficit.
The only letters wich can have a greater deficit are those of $R$,
with $\underline{\lambda}\left(a\right)=\lambda\left(a\right)-c$
so $\lambda\left(a\right)>c$ and are in number $k$ such that $\sum_{a\in R}\lambda\left(a\right)<1$
and then $kc<1$. If one letter of $R$ passes his deficit decreases
and the deficits of the others letters of $R$ increase (see the proof
of lemma $5$); when a letter pass his defect decreases of $\nu$
and the sum of the other deficits increases but increases less than
$\nu$ (because all the excesses decrease and the lemma $1$ ). As
long as this phenomen lasts, no letter acquires excess and all letters
in excess see their excess diminish; as the sum of excesses is equal
to the sum of deficits, the sum of\textbf{ }deficits decreases. No
letter can see his deficit decrease much because this cause other
letters to pass and this deficit decreases a little each time. The
existence of such $N_{i}$ is impossible and no letter can verify
$\underline{\lambda}(a)=\lambda\left(a\right)-c$ if $c$ is the maximum
conceivable value. But in fact it means that no letter cans verify
$\underline{\lambda}(a)<\lambda\left(a\right)$. The theorem is proved.

\medskip{}

\section{Examples and applications.}

\subsection{Approximation of measures on the torus $\mathbb{R}/\mathbb{Z}$ by
measures associated with sequences $(\beta^{n})_{n\geq1}$ modulo
one.}

Here we shall use only finite alphabets in theorem 1, but the consequences
(the very lardge variety of measures on the torus \textbf{$\mathbb{R}/\mathbb{Z}$
}associated with the exponential sequences $\left\{ \beta^{n}\right\} _{n\geq1}$
when $\beta$ travels$\left]1,\infty\right[$) is interesting.

Given a real number $x,$ let $\left\{ x\right\} $ denotes his fractional
part and $\delta_{y}$ denotes the Dirac measure at point $y$.

A sequence $\left(u_{n}\right)_{n\geq1}$ of reals numbers is said
to be distributed according to the probability measure $\nu$ on the
torus if the sequence of measures $\nu_{N}=\frac{1}{N}\left(\delta_{\left\{ u_{1}\right\} }+...+\delta_{\left\{ u_{N}\right\} }\right)$
weakly converge toward $\nu.$ In general a sequences does not admit
a distribution measure but there always exists a set of measures that
we call be associated with this sequence: \textit{we say that the
measure $\gamma$ is associated with the sequence $\left(u_{n}\right)_{n\geq1}$
if it exists an infinite subset $B$ of $\mathbb{N}$ such that $\gamma=\lim_{\begin{array}{c}
N\rightarrow\infty\\
N\in B
\end{array}}\frac{1}{N}\left(\delta_{\left\{ u_{1}\right\} }+...+\delta_{\left\{ u_{n}\right\} }\right).$} 

\textbf{Théorem 2. }\textit{Let $\mu$ be a probability measure on
the torus. then there exists a sequences of real numbers }$\left\{ \beta_{k}\right\} _{k\geq1}$\textit{and
for all $k$ a measure $\mu_{k}$ associated with the sequence $\left(\left\{ \beta_{k}^{n}\right\} \right)_{n\geq1}$
such that $\mu$ is the weak limit of the sequence of measures $\mu_{k}$. }

\textit{In other words the set of associated measure of sequences
$\left(\left\{ \beta^{n}\right\} \right)_{n\geq1}$is dense in the
set of probabilities measures on the torus (in the sense of weak convergence).}

Remark: For almost all $\beta>1$ the sequence $(\beta)_{n>0}^{n}$
is uniformly distributed modulo one. The incriminated \textit{$\beta_{k}$
}are generally part of the exceptional set. 

\textit{Proof} : We will use the following idea:

\textbf{Proposition 4} $\left[1\right]$ \textit{Given a sequence
$\left(I_{m}\right)_{m\geq1}$ of closed intervals $I_{m}$ of the
of the torus $\mathbb{R}/\mathbb{Z}$ , of the same length (or of
length gteater than a positive number $c$), there exists a real number
$\beta$ such that for all $m\geq1$, $\beta^{m}$ belongs to $I_{m}.$}

\textit{Proof of the Theorem.} Let us fix a measure $\mu$ on the
torus; if she has atoms let $(y_{n})_{n\geq1}$ be these atoms; they
are finite or countably in number.

Fix an integer $k,$then let us divide the torus into a finite numbers
of intervals $I_{1}^{k},I_{2}^{k},...,I_{p_{k}}^{k}$whose union is
the entire torus and whose possible intersections two by two are reduced
at their ends. We choose their lengths greater than a number $c_{k}>0$
and such that no atom $y_{i}$, if any, is one of the ends of the
intervals. For any integer $i$ less than or equal to $p_{k}$ we
set $\lambda^{k}\left(i\right)=\mu\left(I_{i}^{k}\right).$ The choice
of the $(I_{i}^{k})_{i=1,...,p_{k}}$ guarantees that the sum of the
$\lambda^{k}\left(i\right)$ is equal to$1$. 

Still for this fixed $k$ apply the theorem $1$ to the alphabet $I_{1}^{k},I_{2}^{k},...,I_{p_{k}}^{k}$,
assigning to $I_{i}^{k}$ the frequence $\lambda^{k}\left(i\right)=\mu\left(I_{i}^{k}\right)$;
construct with the Theorem $1$ a sequence of intervals $(u_{n}^{k})_{n\geq1}$where
$u_{n}^{k}$ is an interval belonging to the alphabet $I_{1}^{k},I_{2}^{k},...I_{p_{k}}^{k}$;$u_{n}^{k}$
take the value $I_{i}^{k}$ with the frequence $\lambda^{k}\left(i\right)$.

According to the proposition $4$ there exists a real number $\beta_{k}$
such that, for all $n,$ $\left\{ \beta_{k}^{n}\right\} $ falls in
the interval $u_{n}^{k}$. The sequence $(\left\{ \beta_{k}^{n}\right\} )_{n\geqslant1}$
does not always admit a distribution measure but she always has an
associated measure $\nu_{k}$; as $(u_{n})=I_{i}^{k}$ with the frequency
$\lambda^{k}\left(i\right)=\mu\left(i\right)$, for all interval $I_{1}^{k},I_{2}^{k},...,I_{p_{k}}^{k}$,
$\nu_{k}\left(I_{i}^{k}\right)=\mu\left(I_{i}^{k}\right)$ (recall
that the intervals $I_{1}^{k},I_{2}^{k},...I_{p_{k}}^{k}$ are disjoints
relatively to the measure $\mu$).

Let us increase $k$ and choose the $I_{1}^{k},I_{2}^{k},...I_{p_{k}}^{k}$
so that the limit of the $\sigma-$algebra generated by the set $\left\{ I_{1}^{k},I_{2}^{k},...,I_{p_{k}}^{k}\right\} $give
rise to the Borelians of the torus, it then comes that $\lim_{k\rightarrow\infty}\nu_{k}=\mu$.

\subsection{Numeration.}

Let $\beta$ be an integer; we shall take $\beta=10$ but $\beta$
could be any integer. Let $(u_{n})_{n\geq1}$ be a sequence wich begins
with a Joker $u_{1}=J$ ($J$ is a symbol wich is only encoutered
there, it correspond to zero: no digit has yet grow but zero is part
of the sequence of numeration in base $10).$ Then we look at the
sequence of expansions of the integers in base $10$ and we write
$u_{N}=a_{i}$ if to go from $N-1$ to $N$ we increase of $1$ the
$i-th$ digit starting from the right (that of the $10^{i-1}$) or
we make it appears if there was none $10^{i-1}$ before (and all the
digits on the right of the $(i-1)th$ digit goes to zero).The beginning
of the sequence is

$u_{1}u_{2}u_{3}u_{4}u_{5}u_{6}u_{7}u_{8}u_{9}u_{10}u_{11}u_{12}u_{13}u_{14}u_{15}u_{16}u_{17}u_{18}u_{19}u_{20}u_{21}u_{22}...u_{101}...u_{111}u_{112}u_{113}...$$=Ja_{1}a_{1}a_{1}a_{1}a_{1}a_{1}a_{1}a_{1}a_{1}a_{2}a_{1}a_{1}a_{1}a_{1}a_{1}a_{1}a_{1}a_{1}a_{1}a_{2}a_{1}...a_{3}a_{1}a_{1}a_{1}a_{1}a_{1}a_{1}a_{1}a_{1}a_{1}a_{2}a_{1}a_{1}...$ 

We obtain the same sequence by constructing the democratic sequence
starting with a Joker (step $1)$ and such that the democratic process
begin at step $2$ with the alphabet $\left(a_{n}\right)_{n\geqslant1}$with
frequency $\frac{9}{10}$ for the symbol $a_{1}$, $(\frac{9}{10})\frac{1}{10}$
for $a_{2}$, $(\frac{9}{10})\frac{1}{10^{2}}$ for $a_{3}$,... $(\frac{9}{10})\frac{1}{10^{n-1}}$for
$a_{n}$ and so on; the $n$- th digit appears $10$ time less than
the$(n-1)$th.

This is also an automatic sequence: let $(a_{n})n\geq1$ an alphabet,
let $\sigma$ be the morphism $a_{1}\rightarrow a_{1}...a_{1}a_{2}$
, $a_{2}\rightarrow a_{1}...a_{1}a_{3},$ and more generally $a_{n}\rightarrow a_{1}...a_{1}a_{n+1}$
(nine $a_{1}$ at each time). Let $\left(u_{n}\right)_{\geq1}$ be
the fix point of the substitution $\sigma$ starting with $a_{1}$
and let $a_{i}$ represent ``the $i$- th digit is growing at step
$n$ if and only if $u_{n}=a_{i}$. The democratic sequence is equal
to the sequence beginning with the Joker followed by the fixed point.

\textbf{Proposition 4.} Let $A=(a_{1},a_{2}....a_{n}...)$ an enumerable
alphabet, and let $\frac{1}{k}(1-\frac{1}{k}),$$\frac{1}{k^{2}}(1-\frac{1}{k})$,...$\frac{1}{k^{n}}(1-\frac{1}{k})$,...a
sequence of numbers with sum $1.$ Let $\left(u_{n}\right)_{n\geq1}$be
the associated democratic sequence beginning by a Joker. If $u_{n}=i,$
the $i-th$ digit is growing at the step $n.$

You can easily check it in the case where $\beta$ is an integer;
we conjecture it for all $\beta>1$ (see $\left[2\right]$ for numeration
in non integer base). 

In the case where this would be true for all $\beta>1,$ maybe it
would help to prove another conjecture: does the $\beta-$shift is
the symbolic system with entropy $log\beta$ having the minimum of
$n-$words? (this is true for $\beta\in\mathbb{N}$ ).

There is also an inteesting problem : is the convergence toward the
right frequencies fast? That might have an impact on the previous
question.

\textbf{Bibliographie}

$\left[1\right]$A. Bertrand Mathis. Sur les parties fractionnaires
des suites $(\beta^{n})_{n\geq1}$, \textit{Uniform Distribution Theory
}\textbf{\textit{14}}\textit{ n\textdegree 2 (2019),69 -72.}

$\left[2\right]$A. Bertrand-Mathis. Comment écrire les nombres entiers
dans une base qui n'est pas entière.\textit{ Acta Math.Hungar. }\textbf{\textit{54}}\textit{(
1989) n\textdegree 3-4 237-241.}

$\left[3\right]$J.A. Ville. Etude critique de la notion de collectif,
Thèse, Faculté des Sciences de Paris, 1939.
\end{document}